\documentclass{amsart}
\usepackage{amssymb}
\usepackage{amscd}
\newtheorem{theorem}{Theorem}[section]

\newtheorem{propo}[theorem]{Proposition}
\newtheorem{corol}[theorem]{Corollary}

\newtheorem{remark}[theorem]{Remark}

\DeclareMathOperator{\sgn}{sgn}
\DeclareMathOperator*{\ooplus}{\oplus}

\DeclareMathOperator*{\too}{\to}
\DeclareMathOperator*{\ssim}{\sim}
\DeclareMathOperator*{\ttimes}{\times}
\DeclareMathOperator*{\ppartial}{\partial}
\DeclareMathOperator*{\Exp}{Exp}
\begin{document}

\title[Universal variations of Hodge structure]
{Universal variations of Hodge structure and Calabi-Yau-Schottky relations}
\author{Ziv Ran}
\address{}
\email{}

\maketitle

\setcounter{section}{0}

Varieties generally come to life through the maps between them and the modules
(e.g. functions) that live on them. Moduli spaces are no exception to this rule.
An important class of maps on moduli spaces (for compact K$\ddot{a}$hler manifolds) is
that of period maps, which are substantially the same thing as the modules known
as variations of Hodge structure.

The purpose of this paper is twofold. First, we give a canonical formula for the
variation of Hodge structure associated to the m-th order universal deformation
of an arbitrary compact K$\ddot{a}$hler manifold without vector fields. Second, we specialize
to the case of a Calabi-Yau manifold $X$ where we give a formula for the m-th
differential of the period map of $X$ and deduce formal defining equations for its
image (Schottky relations); these are (necessarily infinite, in dimension $\geq 3$)
power series in the middle cohomology.

We will use the method of canonical infinitesimal deformations, developed by the
author in earlier papers [R1, R2]. This method gives a canonical description of
infinitesimal moduli spaces and, what's more, natural maps involving them.
While it might be argued that a germ of a smooth space-such as the moduli of an
unobstructed manifold-is a rather rigid featureless object, making a canonical
description of it uninteresting, on the contrary $maps$ involving such germs can
be quite interesting;  in the case of moduli, the method of canonical infinitesimal
deformations provides a vehicle for studying such maps. For instance in the case
at hand the $m$-th derivative of the period map of Calabi-Yau $n$-fold $X$ is
a filtered map
$$T^m_XM \to H^n_{DR}(X)/F^n$$
whose associated gradeds $S^iH^1(T_X) \to H^{n-i,i}_X$ are the so-called Yukawa-
Green forms (cf. [G]). We will develop cohomological formulas for this and other
derivative maps (Theorems
4.1, 4.3  below), which will allow us to determine their image and
derive (Schottky) relations defining this image, essentially in terms of some
generalized Yukawa-Green type forms
(Theorem 4.1, Corollary 4.4) . For $n=2$ we recover the celebrated 'period
quadric' of $K3$ surface theory; for $n \geq 3$ the relations seem to be new.
For $n=3$ this situation is particularly interesting assuming the 'mirror conjecture'
because then the higher derivatives of the period map, here computed, are related
-in fact, carry equivalent information to- the 'quantum cohomology'(esp. numbers
of rational curves, etc.) of the mirror of $X$. We hope to return to this in
greater detail elsewhere.

The present methods should be applicable in other Schottky-type problems: the
case of curves is being developed by G. Liu ( UCR dissertation).

This paper is a revised version of a manuscript entitled 'linear structure on
Calabi-Yau moduli spaces' (May 1993). We are grateful to Professors P. Deligne
and M. Green for their enlightening comments. We are especially grateful to the referee,
for his dedication above and beyond the call of duty, and for his detailed and
insightful comments which have greatly improved this paper.

\section{Modular coalgebra}
In [R2] we gave a characterization of the dual vector space of an artin local
${\mathbb C}$-algebra as a so-called OS (Order-Symbolic) structure. This characterization
was used in constructing the base ring of the universal formal deformation (e.g.
of a compact complex manifold). Now as might be expected, an important role
in the geometry of moduli spaces, both local and global, is played by certain
canonical 'modular' sheaves and modules over them which are naturally associated
to the deformation problem. Thus motivated we now extend the aforementioned dual
characterization from algebras to modules. The latter give rise on the other, '
coalgebra' side to a so-called Modular Order-Symbolic or MOS structure. An
important guiding principle to keep in mind is that the process of passing from
a module $E$ to the associated MOS structure $B(E)$ should involve {\it no
dualising}, i.e. be covariant functorial in E, because we shall want to apply
this when $E$ is a cohomology group and cohomology and dualising don't mix.
Thus while the algebra-OSS correspondence is essentially a matter of brute-force
arrow reversal, the same is not true for module-MOSS.

Now let $(R_n,m_n)$ be an artin local ${\mathbb C}$-algebra of exponent $n$ (which
to us means $m_n^{n+1}=0\neq m_n^n$) and
$V=V^n = m_n^*$ the associated ( standard) OS structure. It will be convenient also
to consider
$$V_0=V^n_0 = R^*_n = V \oplus V^0_0$$
where $V^0_0 = {\mathbb C}1^*$ and $V^i_0 = V^i\oplus V^0_0, i\geq 1$ ( which might be
called an augmented OSS). Given these data, a $V$-compatible MOS structure $B$
consists by definition of a filtered vector space
$$B^0\subset \cdots\cdots\subset B^n=B ,$$
together with a set of symbol maps
$$\sigma^i_B: B^i \to B^i/B^0 \to V^i\otimes_{\mathbb C}B^{i-1}$$
satisfying the obvious 'comodule' rule, which amounts to commutativity of
$$
\begin {CD}
B^i/B^0 @>{\sigma^i_B}>>V^i\otimes B^{i-1} @>{id\otimes \sigma^{i-1}_B}>> {V^i\otimes V^{i-1}\otimes B^{i-2}}\\
@. @V{\sigma^i\otimes id}VV                         @VVV\\
@.  S^2(V^{i-1})\otimes B^{i-1} @>>> V^i\otimes V^{i-1}\otimes B^{i-1},
\end{CD}
$$
$\sigma^i =$ symbol map for $V^i$. A {\it morphism} of MOS structures is
defined in the obvious way.

{\it Example}  Let E be an $R_n$-module and set
$$B^i(E)=B^i_{R_n}(E)=V^i_0\otimes_{R_n}E,\;i=1,\cdots\,n,\; B(E)=B^n(E)$$
(note indeed that $V^i_0$, besides being an OS structure, is also an $R_n$-module
by the obvious rule $(r.f)(x)= f(rx)$),
with symbol map induced from the comultiplication map
$$V_0^i\to V^i \to V^i\otimes V^i_0,$$
dual to $m\otimes_{\mathbb C} R \to m$, by tensoring with $E$. The assignment
$$E \longmapsto B(E)$$
is a covariant functor from R-modules to $V$-compatible MOS structures.

Now it is important to be able to go the other way. To this end define, for an
MOS structure B,
$$C(B) = C^n(B) = {\text Hom}_{MOS}(V^n_0, B).$$
It is easy to check that the $R_n$-module structure on $V_0$ induces one on $C(B)$,
and that $C$ yields a functor from MOS structures to $R$-modules.

We summarise some relevant properties of these functors as follows.

\begin{propo}
(i) $B$ is right exact and $C$ is left exact ;\\
(ii) there are natural maps
$$E \to C\cdot B(E), \hskip 20pt B\cdot C(D)\to D,$$
which are isomorphisms whenever $E$ is free (resp. $D$ is 'cofree', i.e. a sum of
copies of $V_0$.
\end{propo}
\begin{proof}
(i) is clear from the usual exactness properties of Hom and $\otimes$. As for
(ii), the maps are defined by
$$e \longmapsto (v \longmapsto v\otimes e) \in Hom_{MOS}(V_0,V_0\otimes E),$$
$$v\otimes \phi \longmapsto \phi(v),\;\phi \in Hom_{MOS}(V_0, D).$$
It is easy to check that these are well-defined etc. To complete the proof it
suffices to prove that
$$ Hom_{MOS}(V_0, V_0) \simeq R.$$
more precisely that the map $R\longmapsto Hom_{MOS}(V_0, V_0)$,$r\longmapsto rI$,
is a isomorphism. It is easy to see that this map is injective. For surjectivity
we argue by induction. Take $\phi\in Hom_{MOS}(V_0,V_0)$. By induction there exists
$r \in R_n$ inducing the same map as $\phi$ on $V^{n-1}_0$; using compatibility
with comultiplication it is easy to see that $r$ and $\phi$ also induce the same
map on $V=V_0/V^0_0$, and consequently $\phi - rI$ is effectively a map $V^n_0/V^{n-1}_0
\to V^0_0$, i.e. given by an element of $Hom_{\mathbb C}({m_n^n}^*, {\mathbb C}) = m^n_n$,
say $s$. Then $\phi = (r+s)I$.
\end{proof}
{\it Remark 1.2} Predictably the category of $V-$ MOS also admits internal Hom
and tensor product. For example if $B_1,B_2$ are $V-$ MOS we may define an
MOS $C=B_1\otimes_VB_2$ inductively by $C^0=B_1^0\otimes B_2^0$ (all tensor products
over ${\Bbb C}$ unless otherwise mentioned) and setting $C^i$ to be the preimage
of $V^i\otimes C^{i-1}$ by the natural map
$$B_1^i\otimes B_2^i\to V^i\otimes V^i\otimes B^{i-1}\otimes B^{i-1}$$
(where inductively $V^i\otimes C^{i-1}$ sits naturally in $V^i\otimes V^i\otimes
B^{i-1}\otimes B^{i-1} $ via (comultiplication)$\otimes$(inclusion)).
Similarly we may define $sym^j_V(B)$ and $sym^._V(B)$, the 'symmetric coalgebra'
(or OSS) on $B$, which naturally admits an OS structure coming from the natural
maps
$$sym^i_V(B)\to \sum_j sym^j_V(B)\otimes_Vsym^{i-j}_V(B).$$

\section{Subset spaces and some complexes on them}
Recall that , for a topological space $X$, the very symmetric product $X<m>$
introduced in [R1] parametrises nonempty subsets of $X$ of cardinality $\leq m$.
We now introduce an analogous space parametrising subsets with a distinguished
sub-subset. Define
$$X<m,i> = \{(S,T): T\subset S\} \subset X<m+i>\times X<i> ,$$
with the induced topology as (closed) subset; thus $X<m,i>$ is just the graph of the
tautological or incidence correspondence between $X<m+i>$ and $X<i>$. Note
the natural continuous surjective map
$$\pi_{m,i}: X<m>\times X<i> \to X<m,i>,$$
$$(S^\prime, S^{\prime\prime}) \longmapsto (S^\prime\cup S^{\prime\prime}, S^{\prime\prime}).$$
For $i=1$ is easy to see that via this map $X<m,1>$ may  be identified topologically
with the quotient of $X<m>\times X$ by the relation identifying $(S^\prime,x)\in
X<m-1>\times X\subset X<m>\times X$ with $\pi_{m,1}(S^\prime,x)=(S^\prime\cup \{x\},x)
\in X<m-1,1> \subset X<m>\times X.$

For sheaves $A,B$ on $X$-say of modules over some ring which
will typically be $\Bbb C$-we (abusively) denote by
$\lambda^iA\boxtimes B$, as sheaf on $X<i,1>$, the direct image
$\pi_{i,1*}(\lambda^iA\boxtimes B)$.

Now let $g$ be a sheaf of Lie algebras on $X$. In [R1] we associated to $g$ a Jacobi
complex $J_m(g)$ on $X<m>$, an OS structure on $V^m(g) = {\mathbb H}^0(J_m(g))$,
and consequently a (commutative associative) artin local algebra structure on
$R_m(g) = {\mathbb C}\oplus {\mathbb H}^0(J_m(g))^*$. Note that this construction carries over
essentially verbatim to the case where $g$ is a differential graded Lie algebra
(DGLA) sheaf ( i.e. a 'Lie object' is the category of complexes of sheaves on $X$-as
opposed to ordinary sheaves); in the DGLA case $J_m$ becomes a double complex,
but we shall generally identify it and other multiple complexes with the associated
simple complex. See the Appendix for an interpretation of $R_m(g)$ in this case.

Now let $E$ be a $g$-module. We shall associate to $E$ a complex on $X<m,1>$, called
the modular Jacobi complex and denoted by $J_m(g, E)$. The terms are defined by
$ J^i_m(g, E) = \lambda^i(g)\boxtimes E$ on $X<i,1>\subset X<m,1>,\, 0\leq i\leq m$
(where natually $X<0,1>$ is the diagonal $X\subset X\times X = X<1>\times X$),
and the differential $\partial_i: \lambda^i(g)\boxtimes E\to \lambda^{i
-1}(g)\boxtimes E$
is given by the standard formula from Lie algebra homology:
\begin{eqnarray*}
\partial_i(v_1\times\cdots\times v_i\times e) &=&
\frac{1}{i!}\sum_{\sigma\in S_i}(\sgn\sigma)
[v_{\sigma(j)}, v_{\sigma(k)}]v_{\sigma(1)}\times \cdots\times
\hat{v}_{\sigma(j)}\times\cdots\times\hat{v}_{\sigma(k)}
\times\cdots\times v_{\sigma(i)}\times e\\
& &+\frac{1}{i}\sum_{j=1}^i (-1)^j v_1\times\cdots\times
\hat{v}_j\times\cdots\times v_i\times v_j(e).
\end{eqnarray*}
Now we have a natural map
$$\phi: X<m>\times X<m-1,1> \to X<2m,1>$$
$$(S, (S^\prime, x)) \longmapsto (S\cup S^\prime, x)$$
and an evident map of complexes
$$(J_m/J_1)(g,E) \to \phi_*(J_m(g)\boxtimes J_{m-1}(g,E))$$
given by the natural sheaf maps, analogous to the map
$\wedge ^i(V)\to \oplus (\wedge^{i-k}(V)\otimes \wedge^k(V))$ for a vector
space $V$
$$\lambda^i(g)\boxtimes E \to \ooplus_{k=0}^i
\lambda^{i-k}(g)\boxtimes (\lambda^k(g)\boxtimes E)$$
( it is easy to check that this is compatible with differentials). Now let
$p: X<m,1> \to X$ be the natural map. Then by the above the sheaf
$V^m(E):={\mathbb R}^0 p_*(J_m(g,E))$
on $X$ is endowed with an MOS structure with respect to the (constant) OS structure
${\mathbb H}^0(J_m(g))$, whence a sheaf of $R_m(g)$-modules
$C^m(V^m(E))$ which we denote by
$M_m(g,E)$. Note that the assignment $E \to M_m(g,E)$ is a covariant functor from
$g$-modules to $R_m(g)$-modules.

\section{Universal variation of Hodge structure}
Let X be a compact complex manifold with tangent sheaf $T$ so that $H^0(T) = 0$.
In [R1] we constructed the universal formal deformation
$$\hat{X}^u/\hat{R}^u = \lim_m X_m^u/R^u_m\,,$$
where $R^u_m = R_m(T)$ is the algebra associated to the OS structure
${\mathbb H}^0(J_m(T))$. More generally if $g$ is a sheaf of $\Bbb C$-Lie algebra over
a topological space $X$ with $H^0(g)=0$and $E$ is a sheaf of faithful
$g$-modules, we have constructed the universal $m-$th order $g$-deformation $E_m$
of $E$ over $R_m(g)$, the algebra associated to ${\mathbb H}^0(J_m(g))$.
A first point to be made is the the MOS viewpoint permits a more
'conceptual' interpretation of this construction (I am grateful to the referee
for his insistence that this be explained in detail).
\begin{theorem} $E_m$ is canonically isomorphic to $M_m(g,E)$.
\end{theorem}
\begin{proof}
Let $(E^.,\partial)$ be a soft resolution of $E$ and $(g^.,\delta)$ be a
soft resolution of $g$ where $g^.$ is a dgla acting on $E^.$ compatibly with
the $g$-action on $E$. Applying suitable Schur functors as in [FH], $g^.$
induces a soft resolution of $\lambda^i(g)$ which we denote by $(g^._{-i},\delta^._{-i})$,
whence also a soft resolution of $J_m(g,E)$, which may be used to  compute
$p_*(J_m(g,E))$, yielding a complex $(K^.,d^.)$ on $X$ where
$$K^r = \sum_{i,j}\Gamma(g^i_{i+j+r})\otimes E^i, i\geq -m.$$
Now note that, because $H^0(g)=0$, we have $H^j(\lambda^i(g))=0, j<i,$ hence
$K^.$ is acyclic in negative degrees. More precisely, we may 'cancel off'
the negative part of $K^.$ step-by-step as follows. First $d^{-m}$ is a map
$$\Gamma(g^0_{-m})\otimes E^0\to \Gamma(g^1_{-m})\otimes E^0
\oplus \Gamma(g^0_{-m})\otimes E^1\oplus \Gamma(g^0_{-m+1})\otimes E^0 $$
whose first component $\delta^0_{-m}\otimes id_{E^0}$ is already injective.
Consequently $K^.$ is quasi-isomorphic to a complex $K^._{(1)}$ where $$K^{-m}_{(1)}
=0, K^{-m+1}_{(1)}= (\Gamma(g^1_{-m})/B^1(g^._{-m}))\otimes E^0 \oplus \Gamma(g^0_{-m})
\otimes E^1\oplus \Gamma(g^0_{-m+1})\otimes E^0,$$
and $K^i_{(1)}=K^i, i>-m+1.$ Here $B^. $ denotes coboundaries. To be precise
$\Gamma(g^1_{-m})/B^1(g^._{-m})$ is to be thought of as a subspace of
$\Gamma(g^1_{-m})$ complementary to $B^1(g^._{-m})$.
Now the map induced
by $\delta^1_{-m}\otimes id_{E^0}$ on $(\Gamma(g^0_{-m})/B^1(g^._{-m}))\otimes E^0$
is still injective if $m\geq2$ ( again thanks to $H^0(g)=0$), as are $\delta^0_{-m}
\otimes id_{E^1}$ and $\delta^0_{-m+1}\otimes id{E^0}$, hence $K^._{(1)}$ is in turn
quasi-isomorphic to a complex $K^._{(2)}$ in degrees $\geq -m+2$ with
$$K^{-m+2}_{(2)} = (\Gamma (g^2_{-m})/B^2(g^._{-m}))\otimes E^0 \oplus
(\Gamma(g^1_{-m})/B^1(G^._{-m}))\otimes E^1 $$
$$\oplus (\Gamma(g^1_{-m+1})/B^1(g^._{-m+1}))
\otimes E^0\oplus ...$$
Continuing in this manner, we obtain after $m$ steps a
complex $(K^._{(m)},d^._{(m)})$ in nonnegative degrees where $K^0_{(m)}$ is of the form
$\sum_{i\geq 0} G_i\otimes E^i$ for certain vector spaces $G^i$. Now it is easy to
see as above that for $i>0$ the map $\delta\otimes id_{E^i}$ which is a component of
$d^0_{(m)}$ is injective. For $i=0$, $G_0$ will be a certain quotient of
$\sum_0^m sym^j(\Gamma(g^1)$ and among the components of $d^0_{m}$ on $G^0\otimes E^0$,
one lands in  $\Gamma(g^2)\otimes\sum_1^m sym^{j-1}(\Gamma(g^1))\otimes E^0$
(map induced by ($\delta$ plus graded bracket)$\otimes id_{E^0}$).
By definition, the kernel of this
map is precisely $V^m_0\otimes E^0 = ({\mathbb C}\oplus {\mathbb H}^0(J_m(g)))\otimes E^0$
and consequently $K^._{(m)}$ is quasi isomorphic to a complex $K^._{(m+1)}$ in
nonnegative degrees with $K^0_{(m+1)}= V^m_0\otimes E^0$ and, say, $K^1_{(m+1)}=
\sum H^i\otimes E^i$ for certain vector
spaces $H^i$. By similar considerations $K^._{(m+1)}$ is quasi- isomorphic to
a complex $(L^.,\Delta^.)$ which starts
$$ V^m_0\otimes E^0 \to V^m_0\otimes E^1\oplus L^1_0\to ...$$
and where $\Delta^0$ which is induced by $id \otimes\partial^0 +$ (map induced by
$g$-action) clearly goes into $V^m_0\otimes E^1$. Thus ${\mathbb R}^0p_*J_m(g,E)$
is simply given by the kernel of $\Delta^0$. In light of Proposition 1.1 and its
proof, $M_m(g,E)$ coincides with the kernel of a map
$$R_m(g)\otimes E^0\to R_m(g)\otimes E^1$$
given by $id\otimes \partial^0 +$ (map induced by action), which is precisely the
definition of $E_m$.
\end{proof}

For later use we record a corollary of the construction
(which follows easily from the fact that cocycles are locally coboundaries)
\begin{corol} $V^m(E)$ is locally isomorphic to $V^m_0\otimes E$.
\end{corol}

Our purpose now is to extend this construction to the De Rham
complex of $X$ and consequently to obtain, for $X$ K$\ddot{a}$hlerian,
a construction of the universal
variation of Hodge structure associated to $X$, which is a (cohomology) vector
bundle over ${\rm Spec}(\hat{R}^u)$ together with a (Hodge)
filtration and a (
Gauss-Manin) trivialization or flat structure.

Consider the De Rham complex of X
$$\Omega^\cdot: {\mathcal O}_X\too^d \Omega^1_X \to \cdots \to \Omega^n_X .$$
As is well known, $T$ acts on this via Lie derivative $v\times \omega\longmapsto L_v(\omega)$,
an action which commutes with exterior derivative and is a derivation with respect to wedge
products. Consequently, we have an associated modular Jacobi (bi)complex $J_m(T,\Omega^\cdot)$
on $X<m,1>$ which we call the Jacobi-De Rham complex of $X$. This gives rise to a $V^m$-
compatible MOS structure ${\mathbb R}^0{p_*}J_m(T, \Omega^\cdot)$ (i.e a
complex of sheaves of such), whence a complex of $R^u_m = R_m(T)$-modules
$$\Omega^\cdot_m := M_m(T, \Omega^\cdot):M_m(T, {\mathcal O}) = {\mathcal O}_m \to
M_m(T,\Omega^1) \to \cdots .$$
As remarked above, ${\mathcal O}_m$ coincide with ${\mathcal O}_{X_m}$, the structure
sheaf of the $m$-universal deformation $X^u_m/R^u_m$. This identification may be extended
as follows.

\begin{theorem}
In the above situation, $\Omega^\cdot_m$ is canonically isomorphic to the relative De Rham
complex $\Omega^\cdot_{X_m/R_m}$.
\end{theorem}

\begin{proof}
For $\cdot = 0$ we're OK, so consider the case $\cdot = 1$. The $T$-linear derivation
$d: {\mathcal O}\to \Omega^1$ gives rise by functoriality to an $R^u_m$-linear
derivation ${\mathcal O}_m\to \Omega^1_m$, whence an ${\mathcal O}_m$-linear map
$$\phi^1_m : \Omega^1_{X_m/R_m} \to \Omega^1_m\,.$$
Clearly $\phi^1_m\otimes_{R_m}{\mathbb C} = \phi^1_0$ is the identity. Moreover
it is easy to see that both sheaves are locally ${\mathcal O}_m$-free: for $\Omega^1_{X_m/R_m}$
this is obvious and for $\Omega^1_m$ it follows from the fact that
it is $R_m$-flat by construction (indeed a local basis for $\Omega^1$ yields
a map $n{\mathcal O}_m\to \Omega^1_m$ which is clearly surjective and whose
kernel yields 0 when tensored with ${\mathbb C}$, hence is 0 by Nakayama).
Consequently, a
suitable $n\times n$ matrix representing $\phi^1_m$ will have the form $I_n$+
(nilpotent),
hence is invertible, so $\phi^1_m$ is an isomorphism. Now in general, for $\cdot = i$
the derivation property of the Lie derivative action of T
on $\Omega^._X$ makes $\Omega_m^.$ an algebra under wedge products,
whence an ${\mathcal O}_m$-linear
map
$$\phi^i_m = \wedge^i(\phi^1_m): \wedge^i_{{\mathcal O}_m}(\Omega^1_{X_m/R_m}) \to
\Omega^i_m.$$
As both sheaves are locally ${\mathcal O}_m$-free and $\phi^i_o$ is the identity, it follows
as above that $\phi^i_m$ is an isomorphism. Compatibility of $\phi^\cdot_m$ with $d$
-i.e, $d\circ\phi^i_m = \phi^{i+1}_m\circ d$- follows easily from the fact that
this holds for $i = 0$ (by construction) and the multiplicativity of $\phi^\cdot_m$.
This completes the proof.
\end{proof}

Given Theorem 3.2 it is clear in principle that the De Rham cohomology of
$X_m/R_m$ is readable
from the cohomology of the Jacobi-DeRham complex and we shall now make this explicit.
Consider the free $R_m$-module
$$H^r_{DR}(X_m/R_m) = {\mathbb H}^r( X, \Omega^\cdot_{X_m/R_m})\,.$$
This is endowed with a Hodge filtration $F^\cdot$, which is induced by the stupid
filtration on the complex $\Omega^\cdot_{X_m/R_m}$, as well as a (Gauss-Manin) isomorphism
$H^r_{DR}(X_m/R_m) \simeq H^r(X, R_m) \simeq H^r(X, {\mathbb C})\otimes R_m$ induced by
the (quasi-isomorphic) inclusion $R_m \to \Omega^\cdot_{X_m/R_m}$.
When $X$ is K$\ddot{a}$hlerian we have by Deligne [D] that
$F^i/F^{i+1} \cong H^{n-i}(\Omega^i_{X_m/R_m})$, a free $R_m$-module.
The filtered module ($H^r_{DR}(X_m/R_m), F^\cdot$) may be called the $m$-universal
variation of Hodge structure ($m$-UVHS) associated to $X$. Our goal is to give
an explicit formula for it, together with the Gauss-Manin isomorphism, in terms of
$X$ itself.

\begin{theorem}
Let $X_m/R_m$ be the $m$-universal deformation of a compact k$\ddot{a}$hler manifold $X$
with $H^0(T) = 0$. Then:
(i) the Leray spectral sequence
$$E^{p,q}_2 = {\mathbb H}^p(X, R^q{p_*}J_m(T,\Omega^\cdot)) \Rightarrow
{\mathbb H}^{p+q}(X<m,1>, J_m(T, \Omega^\cdot))$$
degenerates at $E_2$;

(ii) we have $R_m$-linear, Hodge filtration-preserving isomorphisms
$$E^{r,0}_2 = E^{r,0}_{\infty} ={\mathbb H}^r(X, B^m\Omega^\cdot_{X_m/R_m}) =
B^m{\mathbb H}^r(X, \Omega^\cdot_{X_m/R_m}) \hskip 20mm (3.1)$$
\begin{eqnarray*}
{\mathbb H}^r(X, \Omega^\cdot_{X_m/R_m}) &\cong &
C^m{\mathbb H}^r(X, {\mathbb R}^0{p_*}R_m(T,\Omega^\cdot)) \hskip 35mm (3.2)\\
& \subset & C^m{\mathbb H}^r(X<m,1>, J_m(T, \Omega^\cdot))\;;
\end{eqnarray*}

(iii) the Jacobi-Hodge-De Rham spectral sequence
$$E^{p,q}_1 = {\mathbb H}^p(X,{\mathbb R}^0{p_*}J_m(T,\Omega^q)) \Rightarrow
{\mathbb H}^{p+q}(X, {\mathbb R}^0{p_*}J_m(T, \Omega^\cdot))$$
degenerates at $E_1$\;;

(iv) the Gauss-Manin isomorphism
$$GM_m: H^r_{DR}(X_m/R_m) \to H^r_{DR}(X)\otimes R_m$$
is adjoint to a map $GM_m: B^mH^r_{DR}(X_m/R_m) \to H^r_{DR}(X)\otimes [{\mathbb C}\oplus
{\mathbb H}^0(J_m(T))]$ induced by a map of complexes
$$M_m: J_m(T, \Omega^\cdot) \to J_m(T,\Omega^\cdot_{triv})\,,$$
$\Omega^\cdot_{triv} = \Omega^\cdot$ with trivial $T$-action, which respects the Hodge
filtration up to a shift of m, i.e. $M_m(F^i) \subset F^{i-m}$, hence $GM_m(F^i)
\subset F^{i-m}$.
\end{theorem}
\begin{proof}
To begin with, the inclusion
$${\mathbb C}_X \to \Omega^\cdot$$
is , by the Poincar\'e lemma, a $T$-linear isomorphism (with the trivial action
of $T$ on ${\mathbb C}_X$). Now recall the finite-to-one map
$$\pi: X<m>\times X \to X<m,1>$$
(cf. section 2) and let $X$ be embedded diagonally $X\hookrightarrow X\times X
= X<1>\times X \subset X<m>\times X$. Then it follows that
$$J_m(T, \Omega^\cdot)\ssim_{qis}J_m(T,{\mathbb C})\ssim_{qis}
 \pi_*(p^*_1J_m(T)\oplus p^*_2{\mathbb C}_X)\,.$$
This clearly implies the degeneration assertion of $(i)$. As for $(ii)$, we have
proven in Theorem 3.1 that
$$\Omega^\cdot_{X_m/R_m} \simeq C^m({\mathbb R}^0p_*J_m(T,\Omega^\cdot)).$$
whence a map
$${\mathbb R}^0p_*J_m(T, \Omega^\cdot) \to B^m\Omega^\cdot_{X_m/R_m}\simeq
B^mC^m{\mathbb R}^0p_*J_m(T,\Omega^\cdot).$$
It is easy to see that this is an isomorphism for $m = 0$, and both sides are
$V^m$-coflat, locally a sum of copies of $V^m_0\otimes {\mathcal O}_X$, hence this
is an isomorphism.
 On the
other hand $\Omega^\cdot_{X_m/R_m}$ being $R_m$-flat we have
$${\mathbb H}^r(X, {\mathbb R}^0p_*J_m(T,\Omega^\cdot)) =
{\mathbb H}^r(X, B^m\Omega^\cdot_{
X_m/R_m}) = B^m{\mathbb H}^r(X, \Omega^\cdot_{X_m/R_m}).$$
This proves (3.1). Applying $C^m$ to both side and using $R_m$-freeness of
${\mathbb H}^r(X,\Omega^\cdot_{X_m/R_m})$, (3.2) follows too. Given $(ii)$,
$(iii)$ follows easily from degeneration of the usual (relative) Hodge-De Rham
spectral sequence
$$E_1^{p,q} = H^q(X,\Omega^p_{X_m/R_m}) \Rightarrow H^{p+q}_{DR}(X_m/R_m).$$
Finally for $(iv)$-which is really the main point-we shall construct an explicit isomorphism based on interior
multiplication. For local vector fields $v_1, \cdots , v_k$ we let
$$i_{v_1\cdots v_k}: \Omega^i \to \Omega^{i-k}$$
denote interior multiplication by $v_1\times \cdots\times v_k$ ( which coincides
with $i_{v_1}\circ \cdots \circ i_{v_k}$). Recall Cartan's formula for the Lie
derivative action of $T$ on $\Omega^\cdot$:
$$L_v(\omega) = i_v(d(\omega)) + d(i_v(\omega))\;.$$
Now let us define a map $M$, preserving total degree
\begin{eqnarray*}
M = M^{\cdots} & :J_m(T,\Omega^\cdot) & \to J_m(T, \Omega^\cdot_{triv}) \hskip 30mm (3.3)\\
M^{i,j,k} & : \lambda^jT\boxtimes \Omega^i & \to \lambda^{j-k}T\boxtimes\Omega^{i-k},
k\geq 0,
\end{eqnarray*}
$$M^{i,j,k}(v_1\times \cdots \times v_j\times \omega)= \pm \sum (-1)^{r_1+\cdots +r_k}
v_1\times\cdots\times\hat{v}_{r_1}\times\cdots\times\hat{v}_{r_k}\times\cdots\times v_j\times
i_{v_{r_1}\cdots v_{r_k}}(\omega)\;.$$
In other words $M$ is the map whose restriction on $\lambda^j T\boxtimes \Omega^i$
is given by $\oplus_{k\geq 0} M^{i,j,k}$; note that each $M^{i,j,0}$ is the 'identity'.
The Cartan formula shows that $ M$ is a morphism of complexes, whence a map
\begin{eqnarray*}
^tGM_m: {\mathbb H}^r(X, {\mathbb R}^0p_*J_m(T,\Omega^\cdot))&\to & {\mathbb H}^r
(X, {\mathbb R}^0p_*J_m(T, \Omega^\cdot_{triv})\\
& = & {\mathbb H}^r(\Omega^\cdot_X) \otimes [{\mathbb C}\oplus {\mathbb H}^0(J_m(T))].
\end{eqnarray*}
In view of the commutative triangle
$$\begin{array}{ccccc}
 & & J_m(T, {\Bbb C}_X) & &  \\
& \wr\swarrow & & \searrow\wr & \\
J_m(T, \Omega^\cdot)& &\stackrel{ M}{\longmapsto}  & & J_m(T,\Omega^\cdot_{triv})
\end{array}
$$
$^tGM_m$ is an isomorphism, and it is obvious that it decreases Hodge level by
at most $m$. Applying the $C^m$ functor now suffices to conclude. Note that by
construction the inverse of our GM isomorphism is induced by the quasi-isomorphism
$R_m\to \Omega^._{X_m/R_m}$; the same is obviously true of the 'usual' Gauss-Manin
in any way it is defined (e.g. [K]). Thus our GM coincides with the usual.

\end{proof}
The following discussion concerning cohomology and obstructions was suggested by
some recent work of Clemens [C]. Consider the exact sequence
$$0\to m^m_m\otimes \Omega^p\to \Omega^p_m \to \Omega^p_{m-1} \to 0. $$
By [D] (assuming $X$ K\"ahlerian) this induces an exact sequence
$$0\to m_m^m\otimes H^q(\Omega ^p)\to H^q(\Omega ^p_m)\to H^q(\Omega^p_{m-1})\to 0,$$
i.e. the coboundary map
$\partial_m :H^q(\Omega ^p_{m-1})\to m^m_m\otimes H^{q+1}(\Omega ^p)$ vanishes, and $H^q(\Omega^p_m)$ is $R_m^u$-free.
It is easy to see that the vanishing of $\partial_{m-1}$ implies a priori that
$\partial_m$ factors through a map $H^q(\Omega^p)\to m_m^m\otimes H^{q+1}(\Omega^p)$.
 'Dually',
let $K^m\subset sym^m H^1(T)$ be the kernel of the $m$-th order
obstruction map $ob_m:sym^m H^1(T)\to {\mathbb H}^1(J_{m-1}(T)$ (cf. [R2],2.3).
Alternatively, $K^m$ may also be defined inductively as the kernel of an obstruction map
$ob_m:K^{m-1}.H^1(T)\to H^2(T)$, where $K^{m-1}.H^1(T)$ denotes the intersection
of $sym^mH^1(T)$ and $K^{m-1}\otimes H^1(T)$ in $\otimes^m H^1(T)$.
Then the natural map $$K^m\otimes H^q(\Omega ^p)\to H^{q+1}(\Omega ^p),$$
which
represents the obstruction to lifting $(p,q)$-cohomology over the universal
$m$-th order deformation, vanishes. Equivalently, the tautological map
$J_m(T,\Omega^p)\to \lambda^m(T)\boxtimes \Omega^p [m]$ induces a map
$V^m(\Omega^p)\to sym^m H^1(T)\otimes \Omega^p$ whose image coincides with
$K^m\otimes \Omega^p$ (cf. Corollary 3.2). Then the induced map on
cohomology
$$H^q(V^m(\Omega^p))\to K^m\otimes H^q(\Omega^p)$$ is surjective, i.e
the coboundary or obstruction map
$$K^m\otimes H^q(\Omega^p)\to H^{q+1}(V^{m-1}(\Omega^p))$$
vanishes (and moreover $H^q(V^m(\Omega^p))$ is cofree). This
follows from [D] together with a trivial 'Universal coefficient
theorem' which says that $H^q(B^m(\Omega^p_m))=B^m(H^q(\Omega^p_m))$
provided $H^{q-1}(\Omega^p_m)$ is free; alternatively one can simply
mimic Deligne's semicontinuity argument for $V^m(\Omega^p)$ in place of
$\Omega^p_m$.

Now recall the calculus fact
(already used in [R3]):
$$i([x,y])(\omega )= L_x(i(y)\omega )-L_y(i(x)\omega )-d(i(x\wedge y)\omega )
+i(x\wedge y)(d\omega ).(*) $$
Define a map $$i:J_m(T)\boxtimes\Omega^p \to J_{m-1}(T,\Omega^{p-1})[1]$$
by the formula
$$(v_1,...,v_k,\omega )\mapsto \sum (-1)^j(v_1,...{\hat v}_j,...,v_k,i(v_j)\omega ).$$
Then (*) shows, assuming $X$ K\"ahler, that $i$ is a 'weak' morphism of
complexes, in the sense that the appropriate diagrams commute in the derived
category. This is good enough to induce a map on cohomology
$$V^m\otimes H^q(\Omega^p)\to H^q(V^{m-1}(\Omega^{p-1}))$$
and similarly for $m-1$, hence a commutative diagram
$$
\begin{CD}
sym^mH^1(T)\otimes H^q(\Omega^p) @>>> sym^{m-1}H^1(T)\otimes H^{q+1}(\Omega^{p-1})\\
@VVV                                @VVV                    \\
{\mathbb H}^1(J_{m-1}(T))\otimes H^q(\Omega^p) @>>> H^{q+2}(V^{m-2}(\Omega^{p-1})\\
\end{CD}
$$
which induces
$$
\begin{CD}
K^{m-1}.H^1(T)\otimes H^q(\Omega^p) @>>> K^{m-1}\otimes H^{q+1}(\Omega^{p-1})\\
@VVV                                   @VVV \\
H^2(T)\otimes H^q(\Omega^p)  @>>>        H^{q+2}(\Omega^{p-1})\\
\end{CD}
$$
As we have seen, the right vertical map vanishes, and we conclude
the following which was first proven by Clemens [C] by another method:
\begin{corol}(Clemens) For $X$ K\"ahlerian, the map
$H^2(T)\to {\rm Hom}(H^q(\Omega^p), H^{q+2}
(\Omega^{p-1}))$ vanishes on the image of the $m$-th order obstruction
map in $H^2(T)$, for all $m$.
\end{corol}

\vskip 2cm
It is instructive to rewrite some of the conclusion of Theorem 3.4 in traditional
geometric language. Thus $H^r_m = H^r_{DR}(X_m/R_m)$ may be viewed as a geometric
vector bundle-viz. $Spec(Sym(H^{r*}_m))$-over $Spec(R_m)=:U_m$ and the Gauss-Manin
isomorphism yields a commutative triangle:
$$\begin{array}{ccccc}
H^r_m & & \simeq & & U_m\times H^r \hskip 3in (3.4) \\
& \searrow & & \swarrow & \\
& & U_m & &
\end{array}
$$
where $H^r = H^r_0 = H^r_{DR}(X)$. The geometric vector bundles $H^r_m$ and
$U_m\times H^r$ correspond to OS structures over $V_m$ which are the 'symmetric
coalgebra' over $V^m$ (see Remark 1.2)
 on the MOS structures $H^r({\mathbb R}^0p_*J_m(T, \Omega^\cdot))$
and $H^r({\mathbb R}^0p_*J_m(T, \Omega^\cdot_{triv}))$, respectively, the
latter being evidently isomorphic to $V^m\otimes S^\cdot(H)$, where $S^\cdot(H)$,
the symmmetric coalgebra (= algebra) on $H$ over
${\mathbb C}$ is endowed with the evident OS structure
(e.g. via its duality with $S^\cdot(H^*))$. Now the Hodge
subbundle $F^i_m \subset H^r_m$, together with the trivialization (3.4), give
rise to a morphism to a Grassmannian
$$p_m: S_m \to {\rm Grass}(\dim F^i_0, H^r)\;,$$
which is none other than the $m$-th order germ  of the period map associated to
$(X, H^r, T^i)$. Thus Theorem 3.2 yields, in principle, an explicit formula for
$p_m$. In the next section we shall make this concrete in the case of Calabi-Yau
manifolds.

\section{The Calabi-Yau case}

In this section we fix a Calabi-Yau manifold $X$, i.e. an $n$-dimensional compact
K$\ddot{a}$hler manifold admitting a nowhere-vanishing $n$-form $\Phi$ (which is then unique
up to a constant), and such that $H^0(T) = 0$. We call such a pair $(X,\Phi)$ a
{\it measured} Calabi-Yau manifold (MCYM). A (local) isomorphism f between MCYM's
$(X,\Phi), (X^\prime, \Phi^\prime)$ is supposed to preserve the $n$-form, i.e.
$f^*\Phi^\prime = \Phi$. Thus the Lie algebra sheaf of infinitesimal automorphisms of
a MCYM $(X,\Phi)$ may be identified as the subset $\hat{T} \subset T$ of divergence-free
vector fields $v$, i.e. those with $L_v(\Phi) = 0 $. The inclusion of Lie algebras
$\hat{T} \subset T$ corresponds to an inclusion of rings
$$R_m = R_m(T) \subset \hat{R}_m = R_m(\hat{T}).$$
Indeed it is easy to see that $\hat{R}_m = R_m[t]/(\underline{m}, t)^{m+1}$. The
universal deformation over $\hat{R}_m$ is just $X^n_m \ttimes_{R_m} \hat{R}_m$.
The advantage of $\hat{R}_m$ is that the pullback of the cohomology bundle $H^n_m
= H^n_{DR}(X^u_m/R_m)$ over $\hat{U}_m = {\rm Spec}(\hat{R}_m)$
admits a 'tautological' section:
namely that corresponding to the map
$$J_m(\Phi): J_m(\hat{T},{\mathbb C}[-n]) \to J_m(\hat{T}, \Omega^\cdot)$$
induced by the $\hat{T}$-linear (!) map ${\mathbb C}[-n] \stackrel{\Phi}{\to} \Omega^\cdot$.
To be precise, there is an inclusion (of a direct summand)
$$\pi_*p^*_1J_m(\hat{T})[-n] \hookrightarrow J_m(\hat{T}, {\mathbb C}[-n])$$
(indeed $\pi_*p^*_1J_m(\hat{T})[-n]$ is identical with the part of $J_m(hat{T},
{\mathbb C}[-n])$ in negative degrees and this part forms a sub, as well as
quotient, complex because the differential in degree -1, given by the action,
vanishes).
This inclusion induces
$$[\Phi]:{\mathbb H}^0(J_m(\hat{T})) \to
{\mathbb H}^n({\mathbb R}^0p_*J_m(\hat{T},\Omega^\cdot))\;,$$
which corresponds to $\Phi$ as a cross-section of the geometric vector bundle corresponding to
$H^n_m\otimes \hat{R}_m$.
On the ring level this corresponds to the $\hat{R}_m$-algebra homomorphism
$$S^._{\hat {R}_m}((H^n_m)^*\otimes \hat {R}_m)\to \hat {R}_m$$
given by the composite (where $S^.=S^._{\mathbb C}$)
$$S^.((H^n_m)^*)\to S^.(\hat{R}_m) \to \hat{R}_m $$
where the first map is $S^.[\Phi]^*$ and the second is $\hat{R}_m$-multiplication.
Now we may follow this by the map induced by the Gauss-
Manin isomorphism $H_m^n \simeq H^n\otimes R_m$ (3.4), then project to the $H^n$ factor.
The map thus obtained
$$\hat{p}_m: \hat{S}_m \to H^n =: H$$
is the $m$-th order period map associated to the MCYM $(X, \Phi)$.
This is related to the usual period  map $p_m$ by the diagram
$$
\begin{CD}
\hat{S}_m @>{\hat{p_m}}>> H\backslash\{0\}     \\
@VVV                       @VVV         \\
S_m @>{p_m}>>          {{\mathbb P}(H)}.
\end{CD}
$$
Now $\hat{p}_m$ corresponds to a homomorphism $\hat{p}^*_m: S^\cdot(H^*) \to\hat{R}_m$,
which is obviously determined by its restriction $\hat{p}^{1*}_m$ on $H^*$, the
linear functions on $H$. Tracing through the construction of the Gauss-Manin isomorphism
in Theorem 3.2, we conclude the following formula for the dual $\hat{p}^1_m$ of $\hat{p}^{1*}_m$:

\begin{theorem}
The map $\hat{p}^1_m: {\mathbb H}^0(J_m(\hat{T}))\to H^n_{DR}(X)$ is given by
${\mathbb H}^0(j_m)$ where $j_m:J_m(\hat{T})\to \Omega^\cdot[n]$ is defined by
$j_m(v_1\times \cdots\times v_k) = i_{v_1\cdots v_k}(\Phi)$.
\end{theorem}
\begin{proof}
Our map is given as the composite of three maps. First, embedding
$J_m(\hat{T})$ in the negative-degree portion of the $n$-th column
(viz. $J_m(\hat{T},\Omega^n) $) of the
double complex $J_m(\hat{T},\Omega^.)$ through multiplication by $\Phi$.
Second, applying our explicit Gauss-Manin (interior multiplication)
operator $M$ constructed in the proof of Theorem 3.2. Third, projecting
to the zeroth row, which is just given by $\Omega^.$. From our formula
for $M$ (3.3), is is plain that this map is given by ${\mathbb H}^0(j_m)$.
\end{proof}

We shall now use this description of the period map to describe a set of defining
equations for its image (Schottky relations). For convenience, set $\hat{T}^m=
{\mathbb H}^0(J_m(\hat{T}))$, $\underline{m}_m = \hat{T}^{m*}$, the maximal ideal
of $\hat{R}_m$. To mirror the adic filtration on $\hat{R}_m$ we define a slight
modification $\hat{F}$ of the Hodge filtration on $H$ (which we identify
with its dual $ H^*$ by Poincar\'e duality) by
\begin{eqnarray*}
\hat{F}^iH =& H  &i=1\\
=&F^iH &i\neq 1.
\end{eqnarray*}
$\hat{F}$ naturally induces a filtration on the symmetric algebra $S(H)$,(=
unique filtration such that the level of $ab =$ (level of $a) +($ level
of $b$)) and we denote
by $S_m$ its $m$-th quotient $S(H)/\hat{F}^{m+1}S(H)$. Note that
e.g. by Theorem 4.1, $\hat{p}^{1*}_m$ on $H$
takes $\hat{F}$ into the adic filtration on $\hat{R}_m$. As $\hat{F}$ was
extended 'multiplicatively' to $S^.(H)$,  $\hat{p}^*_m$ also takes
$\hat{F}$ into the adic filtration, yielding
a map $S_m\to \hat{R}_m$. It is also clear that $\hat{p}^*_m$ induces an isomorphism
$$gr^1_{\hat{F}}S(H) = H/F^2H \simeq gr^1(\hat{R}_m) = \hat{T}^{1*} \;.\hskip 30mm (4.0)$$
It follows firstly that $\hat{p}^*_m$ is surjective('local Torelli') but also,
more significantly that for any $a\in \hat{F}^iH$ we may choose $b\in \ooplus^n_{j=i}S^jH$
so that
$$\hat{p}_n^*(a) = \hat{p}_n^*(b) \in \underline{m}_n^i
\subset \hat{R}_n\;. \hskip 30mm (4.1)$$
Indeed writing $b=\sum_i^n b_j$ we may firstly choose $b_i\in S^i(H)$ by (4.0) so that
$\hat{p}_n^*(a-b_i)\in \underline{m}^{i+1}$, then choose $b_{i+1}\in S^{i+1}(H)$
so that $\hat{p}_n^*(a-b_i-b_{i+1})\in \underline{m}^{i+1}$, etc.
We let $Y_n$ denote the set of elements $a-b$ of this kind as $a$ ranges
over (a basis of ) $F^iH$, $i= 2,\cdots,n$. These are  our
'Schottky relations' of order $n$; they are essentially a refined version
of the well-known Yukawa forms
$$\eta^i: S^i(H^{n-1,1}) \to H^{n-i, i}$$
in that $a$ defines a linear form on $H^{n-i, i}$ and $a\circ\eta^i$ is 'given' by
(the degree-$i$ part of ) $b$ as linear form on $S^i(H^{n-1,1})$.  Next we define
our Schottky relations of higher order. For each
$m\geq n$ we define by induction a lift $Y_m$ of $Y_n$,
$$Y_m \subset I_m = ker(S_m\to \hat{R}_m)\;,$$
as follows. For any $y \in I_m$ let $y^\prime \in S_{m+1}$ be an arbitrary lift of $y$.
Thus $\hat{p}^*_{m+1}(y^\prime) \in \hat{\underline{m}}^{m+1}_{m+1}$ may be expressed
as a polynomial of degree $m+1$ in $\hat{T}^{1*}$, which may be
represented by an element  $z \in S^{m+1}(H)$
and  we clearly have $p^*_{m+1}(z) = p^*_{m+1}(y)$. Now set
$$y^{\prime\prime} = y^{\prime} - z \in I_{m+1}\;\;, $$
$$Y_{m+1}= \{y^{\prime\prime}: y\in Y_m\}\;.$$

\begin{theorem}
$p^*_m$ induces an isomorphism
$$S_m/<Y_m> \simeq \hat{R}_m\; .$$
\end{theorem}
\begin{proof}
We know the induced map is well defined, surjective and takes $\hat{F}$ into the adic
filtration. To prove injectivity it will suffice to show $p^*_m$ induces an
isomorphism, for each i,
$$gr^i(\hat{p}^*_m) : gr^i(S_m/<Y_m>) \to gr^i(\hat{R}) = S^i(\hat{T}^{1*})\;.$$
Given $a\in \hat{F}^iS_m$, it is clear that we may find $b\in S^i(H)$ with
$a-b \in <Y_m>$. But note the exact diagram
$$\begin{CD}
\hat{F}^{i+1}S_m @>>> S^i(H) @>>> S^i(H/F^2H)\\
@. @. @VV{\wr}V\\
@. @. S^i(\hat{T}^{1*})\;.
\end{CD}$$
Thus if $gr^i(\hat{p}^*_m)(a) = 0 $ then $gr^i(\hat{p}^*_m(b)) = 0$. So $b\in \hat{F}^{i+1}(S_m)$.
Hence $a=b=0$ in $gr^i(S_m/<Y_m>)$.
\end{proof}

With more specifically transcendental considerations, we shall next give more explicit
versions of the above results. These involve a certain 'Green-Green' pairing \#,
obtained by combining the Yukawa pairing $*$ with the Green's operator $G$ on $X$.
We now proceed to define \#.

Let $\Phi^{-1}\in H^0(\wedge^nT)$ denote the section dual to the $n$-form $\Phi$, i.e.
$$i_{\Phi^{-1}}(\Phi) = 1\;.$$
For forms $\alpha_i\in A^{n_i}(X)$,$i=1,2$, $\alpha_1\otimes \alpha_2$ may be considered as a tensor
of degree $n_1+n_2$, whence a tensor $i_{\Phi^{-1}}(\alpha_1\otimes \alpha_2)$
of degree $n_1+n_2-n$
and the Yukawa product $\alpha_1 * \alpha_2$ is by definition
the alternation (skew-symmetrization) of $i_{\Phi^{-1}}(\alpha_1\otimes\alpha_2)$. Note that
$$A^{i_1,j_1}*A^{i_2,j_2} \subseteq A^{i_1+i_2-n , j_1+j_2}\;,$$
$$\bar{\partial}(\alpha_1*\alpha_2) = \bar{\partial}\alpha_1*\alpha_2 \pm \alpha_1*
\ppartial^-\alpha_2\; ,$$
and in particular $\alpha_1*\alpha_2$ is $\bar{\partial}$-closed provided $\alpha_1$
and $\alpha_2$ are . Also note that * is defined by local and holomorphic data,
and makes sense for the appropriate sheaves (both holomorphic and $C^\infty$).

Now fixing a K$\ddot{a}$hler metric on $X$ (e.g. the Ricci-flat one), let $G=G_{\bar{\partial}}
=G_{\partial} = \frac{1}{2}G_d$ be the Green's operator on $X$ (cf.  [GH]). Thus
$$(\bar{\partial}\bar{\partial}^* + \bar{\partial}^*\bar{\partial})G(c) = c$$
whenever $c$ is $d-$ or $\partial-$ or $\bar{\partial}$-exact. In fact in
this case we have by the Hodge identities
$$c = \sqrt{-1}\bar{\partial}\partial\Lambda G(c)$$
where $\Lambda$ is the usual 'dual Lefschetz' operator. We define
$$a\#b = \sqrt{-1}\partial\Lambda G\partial (a*b)\;.$$
It is immediate that
$$A^{i_1,j_1}\#A^{i_2,j_2} \subseteq A^{i_1+i_2-n+1,\;j_1+j_2-1}$$
$$\bar{\partial}(a\#b) = \partial(a*b),\; \partial(a\#b) = 0.$$

Now let us identify the tangent sheaf $T$ with $\Omega^{n-1}$ via interior
multiplication, by $\Phi$, and note that the subsheaf $\hat{T}\subset T$ thus corresponds to the
subsheaf of closed forms $\hat{\Omega}^{n-1} \subset \Omega^{n-1}$. The following identity
on $\hat{T}$ is easy to prove but crucial
$$[a,b] = \partial(a*b)\;.\hskip 40mm (4.2)$$
Here $a,b$ can be holomorphic sections of $\hat{T}$ or $\hat{T}$-valued $(0,k)$
forms, i.e. $\partial$-closed $(n-1,k)$-forms. The case $k=1$, i.e. $a,b \in A^
{n-1,1}$, is essentially equivalent to the 'Tian-Todorov Lemma'; note that in
this case
$$a*b = i_{a\wedge b}(\Phi)\;.$$
[We prove (4.2) for $k=0$ as the general case is similar.Write locally
$\Phi = dz_1\wedge ...\wedge dz_n, e_i=\partial/\partial z_i, a=\sum f_ie_i,
b=\sum g_je_j $ where $$\sum \partial f_i/\partial z_i = \sum \partial g_j/\partial z_j
=0  (*)$$
by divergence- freeness. Then
$$a*b = \sum_{i\neq j} (-1)^{i+j}f_ig_jdz_1...\hat{dz_i}...\hat{dz_j}...dz_n,$$
so $\partial (a*b)$ corresponds via $\Phi -$ multiplication to the vector field
(all summations over $i\neq j$)
$$
\sum \partial f_i/\partial z_i g_je_j - \sum \partial f_i/\partial z_j g_je_i
$$
$$-\partial g_j/\partial z_j f_ie_i + \sum \partial g_j/\partial z_i f_ie_j $$
which by (*) equals
$$-\sum_j \partial f_j/\partial z_j g_je_j - \sum_{i\neq j} \partial f_i/\partial z_j
g_je_i $$
$$+\sum_i \partial g_i/\partial z_i f_ie_i - \sum{i\neq j} \partial g_j/\partial z_i
f_ie_j $$
$$=[\sum f_ie_i,\sum g_je_j] = [a,b] $$ \qed ]

Now recall the $m$-th order period map
$$\hat{p}^1_m:\hat{T}^m \to H$$
which we describe in Theorem 4.1 as ${\mathbb H}^0(j_m)$. Note that this description
implies in particular that the $k$-th graded
$$gr^k(\hat{p}^1_m): S^k{\hat{T}^1}\to H^{n-k,k}(X)$$
is induced by projection $\hat{T}^1 = H^{n,0} + H^{n-1,1} \to H^{n-1,1}$ and
Yukawa multiplication $(H^{n-1,1})^k\stackrel{*}{\to} H^{n-k,k}$. We refine this observation
as follows.
\begin{theorem}
There is a natural splitting
$$\hat{T}^m = \ooplus^m_1S^k\hat{T}^1\;,$$
with respect to which the $k$-th component $\hat{p}^{1,k}_m$, i.e. the $k$-th
derivative of the period map, is given by
$$\hat{p}^{1,k}_m(a_1,\cdots,a_k) = \frac{1}{k!}\sum_{\pi\in \sum_k}\sum^k_{j=
max(n-k+1,1)}(-1)^{j-1}(a_{\pi(1)}\#\cdots\#a_{\pi(j)}*a_{\pi(j+1)}*\cdots*a_{\pi(k)})\;(4.3)$$
where the RHS is viewed as a cocycle in $(A^{\cdot\cdot}, \partial, \bar{\partial})$.
\end{theorem}

\begin{proof}
The main thing will be to construct this splitting. The splittings we shall construct for
various $m$ will be mutually compatible, so there is no loss of generality in assuming
$k=m\geq 2$. Let $a_1,...,a_m \in \hat{T}^1 = H^{n,0} + H^{n,1}$. Note
$\Phi$ identifies $H^{n,0}={\mathbb C}$. We shall
construct a lifting $b$ of $a_1\cdots a_m\in S^m(\hat(T)^1)$ to $\hat{T}^m$ .
Let $a^1_i$ be the $(n-1,1)$ component of $a_i$. Then we may write
$a_1\cdots a_m = a_1^1\cdots a_m^1 + a' $ where $a'\in \oplus {j<m} S^i(\hat{T}),$
so the list of $a'$ is defined by induction. Thus we may assume $a_i=a_i^1$
are of type $(n-1,1)$ and harmonic, hence both
$\partial$- and $\bar{\partial}$-closed. Now the cohomology of $J_m(T)$ and $J_m(\hat{T})$
may be computed by formally
applying suitable Schur functors to  Dolbeault complexes resolving $T$
and $\hat{T}$ as in the proof of Theorem 3.1, giving rise in
the case of $T$ to a double complex $(B^{\cdot\cdot}, \delta, \bar{\partial})$, where
$\delta$ is the differential inherited from that of $J_m(T)$. Now we begin with
$$b^{-m,m}= a_1\cdots a_m\in B^{-m,m} = S^m(A^{n-1,1})\;,$$
whose horizontal differential is
$$\delta (b^{-m,m}) = \frac{1}{m!}\sum_{\pi\in \Sigma_m}[a_{\pi(1)},a_{\pi(2)}]a_{\pi(3)}
\cdots a_{\pi(m)}\;.$$
As the $a_i$ are $\partial$-closed, we have
\begin{eqnarray*}
\delta (b^{-m,m}) & = & \frac{1}{m!}\sum_{\pi\in \Sigma_m}\partial (a_{\pi(1)}*a_{\pi(2)})
a_{\pi(3)}\cdots a_{\pi(m)}\\
& = & \frac{1}{m!}\sum_\pi\bar{\partial}(a_{\pi(1)}\#a_{\pi(2)})a_{\pi(3)}\cdots a_{\pi(m)}
\end{eqnarray*}
Thus a natural choice for the component of b in $B^{-m+1, m-1}$ is
$$b^{-m+1, m-1} = \frac{1}{m!}\sum(a_{\pi(1)}\#a_{\pi(2)})a_{\pi(3)}\cdots a_{\pi(m)}\;.$$
With this we indeed have
$$\bar{\partial}(b^{-m+1, m-1}) = \delta(b^{-m,m})\;.$$
Next we must find a $b^{-m+2, m-2}$ with
$$\bar{\partial}(b^{-m+2, m-2}) = \delta (b^{-m+1,m-1})\;.$$
As above , we may set
$$b^{-m+2, m-2} = \frac{1}{m!}\sum_{\pi\in \Sigma_m}((a_{\pi(1)}\#a_{\pi(2)})
a_{\pi(3)})\cdots a_{\pi(m)}\;.$$
Continuing in this way we obtain a hypercocycle $b^{\cdot\cdot}$ for $B^{\cdot\cdot}$,
lifting $a_1\cdots a_m$, yielding the required splitting.

Now by theorem (4.1), it follows that (identifying as always $T$ with $\Omega^{n-1}$):
\begin{eqnarray*}
\hat{p}^1_m(b) & = & \sum^m_{j=1} i_{b^{-m+j, m-j}}(\Phi)\\
 & = & \frac{1}{m!}\sum_{\bar{n}}\sum^{m}_{j=1} i_{(a_{\pi(1)}\#\cdots \#
a_{\pi(j)})\wedge a_{\pi(j+1)}\wedge\cdots\wedge a_{\pi(m)}}(\Phi)\\
 & = & \frac{1}{m!}\sum_{\pi}\sum_{m-j+1 \leq n}(a_{\pi(1)}\#\cdots \#a_{\pi(j)})*
a_{\pi(j+1)}*\cdots *a_{\pi(m)}\;,
\end{eqnarray*}
proving (4.3).
\end{proof}

Using (4.3), we can give a more explicit construction of the Schottky relations
$Y_m$. Let
$$\eta^{i*}= \gamma^i: A^{i,n-i} = A^{n-i,i*}\to S^i(A^{n-1,1})^* = S^i(A^{1,n-1})$$
be the 'Yukawa comultiplication', and likewise
$$\rho^j: A^{1,n-1} \to S^j(A^{1,n-1})$$
be the Green-Green comultiplication, dual to the $\#$ multiplication; we may extend $\rho^j$
as a derivation to a map
$$\rho^{i,j}: S^i(A^{1,n-1}) \to S^{i+j-1}(A^{1,n-1})\;.$$
(i.e. acting on one factor at a time). Then define
$$\nu^{i,j}: A^{i,n-i} \to S^{i+j-1}(H^{1,n-1})$$
as the harmonic projection of $\rho^{i,j}\circ \gamma^i$. (Note that $H^{1,n-1}
=(H^{n-1,1})^*$). It then follows directly from (4.3) and the above recipe for
$Y_m$ that
\begin{corol}
A complete set of $m$-th order Schottky relations for X is given by
$$\{a+\sum^m_{j=1}(-1)^j\nu^{i,j}(a): a\in H^{i,n-i}. i=2,\cdots,n\}$$
\end{corol}

\begin{remark} {\rm Unfortunately the Green's operator $G$ and its relatives, having to
do with explicit realizations of the cohomology, are notoriously difficult to
compute explictly, except in some special cases, e.g. curves. On the other hand
it is possible to work with $\breve{C}$ech cohomology where an explicit analogue
of $G$ has been constructed in the work of Toledo and Tong [TT]. This looks like
a promising way to write down period maps fairly explicitly in general cases.
We hope to return to this elsewhere. See [L] for some explicit calculations in the
curve case.}
\end{remark}

\begin{appendix}
\section{the basic construction, DGLA case}
Lie algebras are important in deformation theory because the symmetries of objects
to be deformed generally form a continuous group whose tangent space (=\{elements infinitely
near identity\}) is a Lie algebra and it is these symmetries that are used in regluing
pieces of the original object to form the deformation. Deformation problems of this
type may be called 'free' or 'unconstrained'. There are however important deformation
problems which are on the contrary 'constrained' or 'semitrivialized' in that
some 'part' of the object to be deformed is to remain undeformed or, more accurately
to deform in a trivialized-rather than just trivial-manner. The most familiar
example of a semitrivialized problem is the Hilbert scheme, i.e. deformations of a submanifold
$Y\subset X$, fixing $X$. The (full) Lie algebra associated to the
embedding $Y\subset X$ is
$T_{Y/X}$, the sheaf of vector fields on $X$ tangent to $Y$ along $Y$, and the corresponding
(unconstrained) deformation problem is that of deformations of the pair $(X,Y)$. The viewpoint
we adopt here is that deforming $Y$ in a fixed $X$ amounts to deforming the pair
$(X,Y)$ and trivializing the $X$ part. In the general case this viewpoint leads to a
(special kind of ) differential graded Lie algebra (DGLA), essentially consisting in a
Lie algebra $g$, a $g$-module $h$ and a derivation $g\to h$; roughly $g$ will do the
deforming and $h$ the trivialising. In the Hilbert scheme case $h=T_X$, of course,
but it is important to observe that the Lie algebra structure of $T_X$ itself plays
no role, only the $T_{X/Y}$-module structure. A succinct way to describe the situation
in this case in terms of the normal sheaf $N=T_X/T_{X/Y}$ is that $N[-1]$ forms a 'Lie algebra
in the derived category' controlling the deformations of $Y$ in the fixed $X$.

Now recall that by definition a $(0,1)$-DGLA consists of  a Lie algebra $g$,
 a
$g$-module $h$ plus a Lie derivation $d: g\to h$; as we shall consider no other kind
of DGLA,
we shall for convenience drop the (0,1) tag. Note that the universal enveloping algebra
$U(g \to h)$ is a differential graded associative algebra, with elements of
$h$ having degree 1 and the differential being the natural extension of $d$ as internal
derivation of degree $+1$, Note that the Lie 'inclusion' $g\to (g\to h)$ induces a
homomorphism $U(g) \to U(g\to h)$ and in view of the rule
$$b.a=a.b - a(b)\;,\;a \in U(g),\;b\in S^\cdot(h)\;,$$
$a(b)$ being the action of $a$ on $b$, it is easy to see that as
$U(g)$-modules, we have
$$U(g\to h) \simeq U(g)\otimes S^\cdot(h).$$

Note that to a DGLA sheaf $g\too^d h$ on $X$ we may associate Jacobi bicomplexes
$J_m(g\to h)$ whose  $j$th row is of the form $\pi_{m-j,j*}J_{m-j}(g, \sigma^jh)$ where
$\pi_{m-j,j}: X<m-j,j> \to X<m>$ is the natural (forgetful) map, zeroth row is
just $J_m(g)$ and whose $k$th column is an Eagon-Northcott type
complex of the form
$$\lambda^kg\to...\lambda^jg\boxtimes \sigma^{k-j}h\to... \sigma^kh;$$
 it is the derivation
property of $d$ that ensures that these fit together
to form a bicomplex. As in the Lie algebra case, ${\mathbb H}^0(J_m(g\to h))$ is
an OS structure with a natural OS morphism ${\mathbb H}^0(J_m(g\to h))\to
{\mathbb H}^0(J_m(g))$,
inducing a ring homomorphism $R_m(g) \to R_m(g\to h)$.
\begin{remark}
{\rm For cohomological purposes we may replace $g\to h$ by a $\breve{C}$ech bicomplex,
thus turning $J_m(g\to h)$ into a
multicomplex $\{\wedge^aC^b(g)\otimes S^cC^d(h)\}$
( where $\wedge$ and $S$ are to be understood in the graded sense). The contribution
to ${\mathbb H}^0$ from the 'pure $h$' part $a=b=0$ may be identified with the
kernel of the differential $\tilde{d}$ induced by $d$ on the enveloping algebra
of the DGLA $C^\cdot(g) \to C^\cdot(h)$. Moreover, in cases of interest to us, we
have that $H^0(g)\to H^0(h)$ is injective and in this case it is easy to check
that there is no other contribution to ${\mathbb H}^0(J)$, so in fact}
$${\mathbb H}^0(J_m(g\to h)) = {\rm coker}  (\tilde{d}). \hskip 30mm (A.1)$$
\end{remark}

Now our purpose is to interpret $R_m(g\to h)$ deformation theoretically. To this end let
$E$ be a $g$-module over $X$, $R_m$ an artin local ${\mathbb C}$-algebra of exponent $m$, and
$E_m$ a flat $g$-deformation of $E$ over $R_m$, i.e. an $R_m$-module locally isomorphic to
$E\otimes R_m$ with gluing maps in $\Exp(g)\otimes R_m = G_m$. As in [R2], $E_m$
comes from a Kodaira-Spencer homomorphism $\alpha_m: R_m(g) \to R_m$, i.e.
$E_m \simeq \alpha^*_mE^u_m$ where $E^u_m$ is the $m$-universal object. Now $h$
itself being a $g$-module, $\alpha_m$ similarly gives rise to a $g$-deformation
$h_m=\alpha_m^*h_m^u$ where
$h_m^u = M_m(g,h)= {\rm Hom}_{MOS}(V^m_0,{\mathbb R}^0p_*(J_m(g,h)$ (Theorem 3.1)(note
objects like $h_m^u$ appear there without the 'u' superscript);
note that $m_mh_m^u = {\rm Hom}_{MOS} (V_m,{\mathbb R}^0p_*(J_m(g,h)$.
Now the canonical map $J_m(g)\to \pi_{m,1*}J_m
(g,h)[1]$ , which is nothing but the map from the zeroth to the first
row of $J_m(g\to h)$, gives rise to a canonical cohomology class $B_m(h)
\in C^m{\mathbb H}^1
(J_m(g,h)) = H^1(h_m^u)$ (in fact this clearly lifts to
$H^1(m_mh_m^u)$ and we shall use $B_m(h)$ to denote
the lift as well), whence a class
$$B(\alpha_m,h) = \alpha_m^*(B_m(h))\in H^1(m_mh_m^u).$$
 Note incidentally that $g$ itself is a $g-$module
via the 1/2 ad action and for this action $g\to^{id}g$ is a
DGLA (the 1/2 factor is needed to make $id$ a derivation),whence
a canonical element $B(\alpha_m,g)\in H^1(m_mg_m).$
For instance if $g= T_X$ then $g_m=T_{X_m/R_m}$,
the relative tangent sheaf.

By a $(g\to h)$-deformation of $E$ we mean the data $(\alpha_m, A)$
consisting of a $g$-deformation of $E$ given by $\alpha_m$ plus a 'trivialization' of
the corresponding canonical cohomology class, given e.g. by a $\check{C}$ech
cochain $A = (v_\alpha) \in \check{C}^\circ(h_m)$, $\delta A= B(\alpha_m,h)$.
These objects can all be represented concretely as in [R2].
Recall to begin with that $\alpha_m$ corresponds to a morphic
hypercocycle $v_m=\epsilon (u_m)\in {\mathbb H}^0(J_m(g))\otimes m_m$
where $u.=u_m\in {\breve C}^1(g)\otimes m_m$ satisfies the ${\breve C}$ech
integrability condition
$$
\delta (u_m) =-1/2 [u_m,u_m];\hskip 20mm (A.2)$$
a deformation corresponding to $\alpha _m$, such as $E_m$ (or
$g_m$ or $h_m$) is essentially {\it defined} by the condition
that its ${\breve C}$ech  complex $({\breve C}^.(E_m),\delta _m)$
be isomorphic to  $({\breve C}^.(E)\otimes R_m,\delta_0\otimes 1
+u_m)$ (viewing $u_m$ as an operator of degree +1 on
${\breve C}^.(E_m)\otimes R_m$). Then $B(\alpha_m,g)$ is represented by
$(u_m)$ (which by (A.2) is indeed a cocycle for $g_m$ and in fact
for $m_mg_m$, which is a flat $R_{m-1}$-module);
similarly, $B(\alpha_m,h)$ is
represented by  $d(u_m).$
Put another way,
the operator $\delta_0\otimes 1$ on
${\check C}^(E)\otimes R_m$ pulled over to ${\check C}^.(E_m)$
becomes $\delta - u_m$, so $u_m$ is indeed a cocycle for $m_mg_m$.
Consequently $A$ may
be represented by
$$ a.\in {\breve C}^0(h)\otimes m_m,$$
$$\delta (a.)+u.(a.) = d(u.)$$
i.e. $a_{\alpha} -a_{\beta} +1/2 u_{\alpha\beta}(
a_{\alpha}+a_{\beta}) = d(u_{\alpha\beta}).$
The following then is the DGLA analogue of the main result
of [R2].
\begin{theorem}
Assume $H^0(g)\to H^0(h)$ is injective. Then for any $m$, there exists a $(g\to h)$
deformation $(\alpha_m^{uc}, A_m^{uc})$ over $R_m^{uc} = R_m(g\to h)$ which is
constrained universal in that for any $(g\to h)$ deformation $(\alpha_m, A_m)$
over $R_m$, there is a factorization of $\alpha_m$
$$R^u_m \too^{\alpha^{uc}_m} R_m^{uc} \too^{\beta_m} R_m$$
such that $A_m = \beta^*_m A^{uc}_m$.
\end{theorem}
\begin{proof}
Analogous, mutatis mutandis, to that in [R2] (compare Theorem 3.1 above).
To begin with, $\alpha^{uc}_m$ simply corresponds to the natural map $q_m :J_m
(g\to h) \to J_m(g)$, while $A^{uc}_m$ comes from the extra data needed to lift
an element of ${\mathbb H}^0(J_m(g))$ to one in ${\mathbb H}^0$ of the double
complex formed from the last 2 rows of $J_m(g\to h)$. Now given an arbitrary
$(g\to h)$-deformation $(\alpha_m, A_m)$ over $R_m$, $\alpha_m$ corresponds to
a cochain $u_m \in C^1(g)\otimes m_m$ and $A_m$ to $a \in C^0(h)
\otimes m_m$ as above. Then
$$(u_m, \frac{1}{2}u^2_m,\cdots,\frac{1}{m!}u^m_m;a, u_m\times a,\cdots, \frac{1}{(m-1)!}
u^{m-1}_m\times a; \frac{1}{2}a^2, u_m\times \frac{1}{2}a^2;\cdots;\cdots;\frac{1}{m!}a^m)$$
$$\hskip 4.4in (A.2)$$
yields a morphic element of $m_m\otimes {\mathbb H}^0(J_m(g\to h))$ corresponding
to $\beta_m$ as required.
\end{proof}

{\it Examples}

1. Let $X$ be a compact complex complex manifold and $Y\subset X$ a submanifold
with normal bundle $N$. We have a DGLA
$$(T_{X/Y} \to T_X) \sim N[-1]\;.$$
As noted before ([R2],section 5), $T_{X/Y}$ controls deformations ($Y_m \subset
X_m/R_m$) of the pair $Y\subset X$ and given such a deformation with Kodaira-Spencer
homomorphism $\alpha_m: R_m(T_{X/Y}) \to R_m$,
evidently $B(\alpha_m,T_X)$ is just the
class corresponding to the induced deformation $X_m$ of $X$ alone,
so that the
data $A$ is simply an identification of $X_m$ with the trivial deformation
$X\times {\rm Spec}(R_m)$.

\begin{remark}
{\rm The cokernel of $\tilde{d}$ on $U(N[-1])$ is a certain quotient of $S^\cdot(T_X)$
which is essentially the sheaf $D^\cdot_{Y\to X}$ of 'normal differential operators'
along $Y$, considered by Burchard, Clemens et al [BC]. Note the natural map}
\begin{eqnarray*}
{\mathbb H}^0(J_m(N[-1])) = {\rm coker}(\tilde{d}, {\mathbb H}^0) & \to & {\mathbb H}^0
({\rm coker}(\tilde{d}, U(N[-1])))\\
&\to & {\mathbb H}^0(D^\cdot_{Y\to X}).
\end{eqnarray*}
\end{remark}

2. (cf.[R3]) Let X be a compact complex K$\ddot{a}$hler manifold with tangent sheaf $T$
and $\eta \in H^{p,q}(X)$. We have a map
$$T \to \Omega^{p-1}[q] \hskip3in (A.3)$$
given by interior multiplication by $\eta$. Representing $\eta$ by a closed form,
note the formula
$$i_{[x,y]}(\eta) = L_x(i_y(\eta)) - L_y(i_x(\eta))-d(i_{x\wedge y}(\eta))$$
for vector fields $x$,$y$. Replacing $\Omega^{p-1}$ by, e.g. its Dolbeault resolution
we have by the usual K$\ddot{a}$hler machinery
$$d(i_{x\wedge y}(\eta)) = \bar{\partial}(j_{x\wedge y}(\eta))$$
for a suitable well-defined $j_{x\wedge y}(\eta) \in A^{p-1,q-1}$, therefore this
term vanishes in the derived category and $
{\mathcal L}_\eta = (T\to \Omega^{p-1}[q])$
forms a DGLA. As shown in [R3], $
{\mathcal L}_\eta$ controls precisely the deformations of
X in which (the GM-constant lift of) $\eta$ maintains Hodge level $p$. In the
situation of example 1, we may take $\eta=[Y] \in H^{p,p}(X)$, $p=
{\rm codim}(Y)$ and
then we have an exact diagram
$$\begin{CD}
N[-1]             @>>> T^\prime     @>>> T @>>>\\
@VVV                      @VVV           @|    \\
\Omega^{p-1}[p-1] @>>> {\mathcal L} @>>> T @>>>.
\end{CD}$$
See [R3] for geometric application of this.

\end{appendix}

\end{document}